\input amstex
\documentstyle{amsppt}
\magnification=1200
\NoBlackBoxes
\hsize=5.5in
\vsize=7.5in
\topmatter
\title One Dimensional Nearest Neighbor Exclusion Processes in Inhomogeneous and Random
Environments\endtitle
\rightheadtext {Inhomogeneous Exclusion Processes}
\date December 21, 2006\enddate
\author Lincoln Chayes and Thomas M. Liggett\endauthor
\thanks Research supported in part by
NSF Grant DMS-0306167 (Chayes) and  DMS-0301795 (Liggett).\endthanks
\affil University of
California, Los Angeles\endaffil
\abstract The processes described in the title always have reversible
stationary distributions. In this paper, we give sufficient conditions
for the existence of, and for the nonexistence of, nonreversible stationary distributions.
In the case of an i.i.d. environment, these combine to give a necessary and sufficient condition for
the existence of nonreversible stationary distributions.
\endabstract 
\endtopmatter

\heading 1. Introduction\endheading

Transport phenomena for noninteracting particles in one dimensional environments is
a well studied subject in the contexts of classical and quantum systems. However,
the influence of interactions among particles in these situations is considerably
less well understood. In this paper, we consider one of the more important models of
particle motion with interaction -- the exclusion process. In particular, we will study
the exclusion process on $Z^1$ with nearest neighbor jumps with probabilities
$p_i\in (0,1)$ and $q_i=1-p_i$ from $i$ to $i+1$ and $i-1$ respectively. This is the continuous
time Markov process $\eta_t$ on $\{0,1\}^{Z^1}$ with formal generator 
$$\Cal L f(\eta)=\sum_i\bigg\{\eta(i)[1-\eta(i+1)]p_i+\eta(i+1)[1-\eta(i)]q_{i+1}\bigg\}
[f(\eta_{i,i+1})-f(\eta)],$$
where $\eta_{i,i+1}$ is obtained from $\eta$ by interchanging the $i$th and $(i+1)$st coordinates. 
The exclusion process has been the object of a lot of attention over the past 35 years,
primarily in case the transition probabilities are translation invariant. 
Here, as the title indicates, we investigate problems where the $p_i$ are inhomogeneous,
with particular applications to the case in which the $p_i$ are i.i.d. random
variables. Among the few  rigorous papers
dealing with spatially inhomogeneous asymmetric exclusion processes are Liggett (1976) and Jung (2003).

The i.i.d. model was investigated in a certain approximate fashion  
in Tripathy and Barma (1998) and, more recently, in Harris and  
Stinchcomb (2004).  In general, these authors have studied the  
``phase diagram" of (maximal) current flow as a function of an  
equilibrium particle density parameter in the presence of disorder.   
Certain interesting phenomena have been uncovered in these works. In  
particular, a symmetric flat region in the above mentioned response  
curves, indicating a  forbidden interval of densities, and the  
observation, primarily numerical, in the earlier reference that in  
the presence of a {\it random} locally preferred direction of flow,  
the current density vanishes with increasing system size.    
Our results bolster some of these conclusions:  We demonstrate the  
existence of current carrying states whenever, for some $\epsilon > 0 
$,  $p_0 \geq \frac{1}{2} + \epsilon$.  Moreover, the flux in these  
states vanishes linearly with $\epsilon$.  Furthermore, we vindicate  
conclusively the later phenomenon.  In particular, whenever $p_0 < q_0 
$ and $p_0>q_0$ have positive probability -- or even if $p_0 = \frac{1}{2}$ is in  
the support of the disorder distribution -- we show that the current  
indeed vanishes  as the system size tends to infinity.  Furthermore,  
under the general condition of zero current, we can characterise all  
the invariant measures.

This picture is in sharp contrast to the non-interacting
version of this problem:
At the end of this paper, we will consider briefly the non-interacting case, and show
that there is a stationary distribution for the system with nonzero flux whenever
the $p_i$'s are i.i.d. and $E\log(p_i/q_i)$ exists and is nonzero. This difference
between the interacting and noninteracting systems might be a bit surprising, since
at least at low densities, one might expect the two systems to have similar
properties.

Returning to the exclusion process, we first note that reversible stationary distributions always
exist in our situation. To define them, let
$\pi_i$ be defined by taking $\pi_0>0$ and then $\pi_ip_i=\pi_{i+1}q_{i+1}$ for all $i\in Z^1$. 
The corresponding reversible measure is the product measure $\nu_{\alpha}$, where
$$\nu_{\alpha}\{\eta:\eta(i)=1\}=\alpha_i=\frac{\pi_i}{1+\pi_i}.\tag 1.1$$
(See Section VIII.2 of Liggett (1985), for example.) According to Theorem 2.1 of Jung (2003),
these are extremal in the class $\Cal I$ of all stationary distributions if and only if
$$\sum_i\alpha_i(1-\alpha_i)=\infty.\tag 1.2$$
When this sum is finite, extremal stationary distributions $\mu_n$ (with
$n\geq 0$ or $-\infty<n<\infty$, depending on the situation) are
constructed from these product measures in the following way: $\mu_n(\cdot)=\nu_{\alpha}(\cdot\mid
A_n),$ where
$$A_n=\cases\{\eta:\sum_i\eta(i)=n\}&\text{ if }\sum_i\alpha_i<\infty,\\
\{\eta:\sum_i(1-\eta(i))=n\}&\text{ if }
\sum_i(1-\alpha_i)<\infty, \\ \{\eta:\sum_{i\in T}\eta(i)-\sum_{i\notin T}(1-\eta(i))=n\}
&\text{ otherwise,}
\endcases$$
where in the third case, $T$ is chosen so that $\sum_{i\in T}\alpha_i<\infty$ and 
$\sum_{i\notin T}(1-\alpha_i)<\infty$.  These conditional measures do not depend on
$\alpha$, since an irreducible positive recurrent Markov chain has a unique stationary
distribution. In the third case, $\{\mu_n\}$ does not depend on the choice of $T$, except for a
possible relabelling. 

In the spatially homogeneous case $p_i\equiv p$, the extremal stationary distributions are
completely known (Liggett (1976)):
$$\Cal I_e=\{\nu_{\rho},0\leq \rho\leq 1\}$$
if $p=1/2$ and 
$$\Cal I_e=\{\nu_{\rho},0\leq \rho\leq 1\}\cup\{\mu_n,-\infty< n<\infty\}$$
if $p\neq 1/2.$ Note that in the latter case, $\nu_\rho$ is not reversible.
(Here $\nu_{\rho}$ denotes the homogeneous product measure of density $\rho$.)

The spatially inhomogeneous case in which $Z^1$ is replaced by $\{0,1,2,...\}$ was
treated in Liggett (1976). In that case, the result is that all stationary
distributions are reversible. 

Our main objective in this paper is to say what we can about the following question:
In the spatially inhomogeneous case on $Z^1$, when do nonreversible stationary
distributions exist? We will give sufficient conditions for the existence and for
the nonexistence of such distributions; they become necessary and sufficient in 
the case of i.i.d. $p_i$'s.

An important tool in discussing this issue (as in many 
involving the exclusion process) is the flux. For $\mu\in \Cal I$, this is defined by
$$\phi(\mu)=p_i\mu\{\eta:\eta(i)=1,\eta(i+1)=0\}-q_{i+1}\mu\{\eta:\eta(i)=0,\eta(i+1)=1\}.\tag 1.3$$
The fact that this quantity is independent of $i$ can be checked by using 
$\int\Cal L fd\mu=0$ for $f(\eta)=\eta(i).$  In case $p_i\equiv
p$, the flux under $\nu_{\rho}$ is $(p-q)\rho(1-\rho)$, for example.

It is easy to check that the flux is zero
for the reversible stationary distributions described above.  Our first result is
a converse to this observation.

\proclaim {Theorem 1} Suppose $\mu\in \Cal I_e$ and $\phi(\mu)=0$. If $\mu$
is not the pointmass on $\eta\equiv 0$ or on $\eta\equiv 1$, then $\mu=\nu_{\alpha}$ 
for some $\alpha$ satisfying
(1.1) if (1.2) is satisfied and $\mu=\mu_n$ for some $n$ if not.\endproclaim

For applications of Theorem 1, it is useful to make the following observation:

\proclaim{Proposition 1} Take $\epsilon>0.$ Suppose that for each $n$ there exists a $k$ so that
$p_i\geq \frac 12-\epsilon$ for all $k\leq i\leq k+n$. Then $\phi(\mu)\geq -\epsilon$ for all $\mu\in
\Cal I$. Similarly, if for each $n$ there exists a $k$ so that
$p_i\leq \frac 12+\epsilon$ for all $k\leq i\leq k+n$, then $\phi(\mu)\leq \epsilon$ for all $\mu\in
\Cal I$.
\endproclaim

Theorem 1 and Proposition 1 will be proved in Section 2. 
Combining these two results, we obtain a sufficient condition for all 
stationary distributions to be reversible:

\proclaim {Corollary 1} Suppose that for every $\epsilon>0$ and every 
positive integer $n$ there exist $k$ and $l$ so
that
$p_i\geq \frac 12-\epsilon$ for all $k\leq i\leq k+n$ and $p_i\leq \frac 12+\epsilon$ for all $l\leq
i\leq l+n$. Then $\mu\in\Cal I$ implies that $\mu$ is reversible.
\endproclaim

This conclusion that all stationary distributions are reversible under the 
much stronger assumption that
$$\lim_{i\rightarrow\pm\infty}p_i=\frac 12$$
is a consequence of Theorem 1.1 of Jung (2003). This conclusion also follows from  Theorem 1.2 in that
paper (assuming $\inf_i p_i>0,\inf_i q_i>0$) in all cases {\it other than} (a)
$\lim_{i\rightarrow -\infty}\pi_i=0, \lim_{i\rightarrow +\infty}\pi_i=\infty$ or
(b) $\lim_{i\rightarrow -\infty}\pi_i=\infty, \lim_{i\rightarrow +\infty}\pi_i=0$.
When the $p_i$'s are i.i.d., these excluded cases are of course the prevalent ones:
case (a) occurs when $E\log (p_0/q_0)>0$ and case (b) when $E\log (q_0/p_0)>0$.

To prepare for the next result, we define the exclusion process on $[m,n]=\{m,...,n\}$ with
boundary conditions by allowing the usual transitions in $\{m,...,n\}$ together with
the transitions $0\rightarrow 1$ at site $m$ at rate $p_{m-1}$ and $1\rightarrow 0$
at site $n$ at rate $p_n$. This is a finite state irreducible Markov chain, and hence
has a unique stationary distribution $\mu_{m,n}$. 

\proclaim{Theorem 2} The flux $\phi(\mu_{m,n})$ is an increasing function of
$p_{m-1},p_m,...,p_n$.\endproclaim

Combining this with the known behavior of $\mu_{m,n}$ in the
homogeneous case yields the following sufficient condition for the
existence of a nonreversible stationary distribution:

\proclaim{Corollary 2} If for some $\epsilon>0$, $p_i\geq \frac 12+\epsilon$ for all $i$, then there
exists a (nonreversible) $\mu\in \Cal I$ with $\phi(\mu)\geq \frac{\epsilon}2$.\endproclaim

This answers a question raised near the end of the introduction to Jung (2003).
Theorem 2 and Corollary 2 will be proved in Section 3. 

Combining Corollaries 1 and 2, we have the following result for the exclusion process
in which the $p_i$'s are chosen randomly in an i.i.d. fashion:
 
\proclaim{Theorem 3} Consider an exclusion process with i.i.d. $p_i$'s. The following
hold with probability 1:

(a) If for every $\epsilon>0$, $P(p_0\geq \frac 12-\epsilon)>0$ and $P(p_0\leq \frac 12+\epsilon)>0$,
then all stationary distributions are reversible.

(b) If for some $\epsilon>0$, $P(p_0\geq\frac 12+\epsilon)=1$ or $P(p_0\leq\frac 12-\epsilon)=1$, then
there exists a nonreversible stationary distribution.
\endproclaim

\noindent {\bf Remark:}
Part (b) clearly follows from Corollary 2 if we replace the i.i.d. assumption
with the assumption that $\{p_i, i\in Z^1\}$ be stationary and ergodic. The same is
not true for part (a). To see this, consider the case of deterministic $p_i$'s with
$$p_i=\cases \alpha&\text{ if } i\text { is even},\\
\beta&\text{ if } i\text { is odd.}\endcases\tag 1.4$$
Then all stationary distributions are reversible if and only if $\alpha+\beta=1$.
Indeed, if $\alpha+\beta=1$, then the process is symmetric, so all stationary 
distributions are exchangeable (and therefore reversible in this case) by Theorem 1.12 of Chapter
VIII of Liggett (1985). On the other hand, the homogeneous product measures $\nu_{\rho}$
are stationary for all choices of $\alpha$ and $\beta$ by Theorem 2.1(a) of the same
chapter. However, since
$$\phi(\nu_{\rho})=(\alpha+\beta-1)\rho(1-\rho),$$
$\nu_{\rho}$ is not reversible if $\alpha+\beta\neq 1$. By letting $p_i$ be
given by (1.4) with probability $\frac 12$ and its translate with probability
$\frac 12$, one obtains a stationary, ergodic sequence with nonreversible
stationary distributions. So, in order to conclude that all stationary
distributions are reversible, it is not enough to assume that $p_i>\frac 12$ and
$p_i<\frac 12$, each a positive proportion of the time. 

\heading 2. Sufficient conditions for all stationary distributions to be
reversible\endheading

In this section, we prove Theorem 1 and Proposition 1. 

\demo{Proof of Theorem 1} The proof is based on the proofs of Theorems 1.2 and 1.3
of Liggett (1976), so we will only sketch the parts that are similar. 
We will use the basic coupling of two copies
$\eta_t$ and $\eta_t'$ of the process. In this
coupling, particles move together as much as possible. At rate $p_i$, for example,
$(\eta,\eta')\rightarrow$
$$\cases(\eta_{i,i+1},\eta')&\text{ if }\eta(i)=1,\eta(i+1)=0;\ 
\eta'(i)=\eta'(i+1)\text{ or }\eta'(i)=0,\eta'(i+1)=1,\\
(\eta,\eta'_{i,i+1})&\text{ if }\eta'(i)=1,\eta'(i+1)=0;\ 
\eta(i)=\eta(i+1)\text{ or }\eta(i)=0,\eta(i+1)=1,\\
(\eta_{i,i+1},\eta'_{i,i+1})&\text{ if }\eta(i)=\eta'(i)=1;\ 
\eta(i+1)=\eta'(i+1)=0.\endcases$$
Quantities related to the coupled process will be denoted by tildes. For example,
the set of stationary distributions for the coupled process will be called $\tilde
\Cal I.$

For $m<n$, define the following two functions of a coupled configuration:
$$\align f_{m,n}(\eta,\eta')=&\sum_{k=m}^n|\eta(k)-\eta'(k)|,\\
g_{m,n}(\eta,\eta')=&\#\text{ of strict sign changes in the sequence }\{\eta'(m)-\eta(m),...,
\eta'(n)-\eta(n)\}.\endalign$$
The fundamental property of these functions that makes them useful is that they cannot
increase except as the result of transitions across the boundaries of $[m,n]$.
Interior transitions can make them decrease unless $\eta_t\leq\eta_t'$ or $\eta_t\geq\eta_t'$
in the case of $f_{m,n}$, and unless there is at most one sign change in the case of
$g_{m,n}$.

We will use the following notation. If $\nu$ is a probability measure on $\{0,1\}^{Z^1}\times
\{0,1\}^{Z^1}$, then $\nu\{(\eta,\eta'):\eta(i)=\delta,\eta'(i)=\delta'\}$ will be denoted by
$$\nu\left(\matrix \delta'\\\delta\\i\endmatrix\right),$$
with analogous notation for probabilities of cylinder sets involving more than one site.

If $\nu\in\tilde\Cal I$, then $\int\tilde\Cal Lf_{m,n}d\nu=0.$ Writing this out gives
$$\gathered 2\sum_{i=m}^{n-1}(p_i+q_{i+1})\bigg[\nu\left(\matrix 1&0\\0&1\\i&i+1\endmatrix\right)+
\nu\left(\matrix 0&1\\1&0\\i&i+1\endmatrix\right)\bigg]=\\
p_{m-1}\bigg[\nu\left(\matrix 1&0\\0&0\\m-1&m\endmatrix\right)-\nu\left(\matrix
1&1\\1&0\\m-1&m\endmatrix\right)+
\nu\left(\matrix 0&0\\1&0\\m-1&m\endmatrix\right)-\nu\left(\matrix
1&0\\1&1\\m-1&m\endmatrix\right)\bigg]\\+ q_{m}\bigg[\nu\left(\matrix 1&1\\0&1\\m-1&m\endmatrix\right)
-\nu\left(\matrix 0&1\\0&0\\m-1&m\endmatrix\right)+
\nu\left(\matrix 0&1\\1&1\\m-1&m\endmatrix\right)-\nu\left(\matrix
0&0\\0&1\\m-1&m\endmatrix\right)\bigg]\\\\ +\text{ similar terms coming from the right boundary of
}[m,n].\endgathered\tag 2.1$$ 
The left side of (2.1) represents interior transitions that lead to the loss of (two) discrepancies,
while the right side corresponds to the gain or loss of discrepancies in $[m,n]$ due to transitions
across the boundary.

At this point, the argument differs somewhat according to whether one,
both, or neither of the following hold:
$$\inf_{i<0}(p_i+q_{i+1})>0,\quad \inf_{i>0}(p_i+q_{i+1})>0.\tag 2.2$$
If both of these infima are 0, then one can let $m\rightarrow-\infty, n\rightarrow
+\infty$ in (2.1) along appropriate subsequences to conclude that 
$$\sum_{i=-\infty}^{\infty}(p_i+q_{i+1})\bigg[\nu\left(\matrix 1&0\\0&1\\i&i+1\endmatrix\right)+
\nu\left(\matrix 0&1\\1&0\\i&i+1\endmatrix\right)\bigg]=0.\tag 2.3$$
It follows from this that $\nu$ puts no mass on configurations with adjacent discrepancies
of opposite types, and then using the invariance of $\nu$ again, that it puts no mass
on configurations with discrepancies of opposite type at all. Therefore, it follows that
$$\nu\{(\eta,\eta'):\eta\leq\eta'\text{ or } \eta'\leq \eta\}=1.\tag 2.4$$

We will now assume that (2.2) holds. (The argument in the case that one of the infima
in (2.2) is zero and the other is positive is a combination of these two arguments, and
will be omitted.) It now follows that
$$\sum_{i=-\infty}^{\infty}\bigg[\nu\left(\matrix 1&0\\0&1\\i&i+1\endmatrix\right)+
\nu\left(\matrix 0&1\\1&0\\i&i+1\endmatrix\right)\bigg]<\infty.$$
Using the invariance of $\nu$ again, it follows that 
$$\sum_{i=-\infty}^{\infty}\bigg[\nu\left(\matrix 1&0\\0&1\\i&i+k\endmatrix\right)+
\nu\left(\matrix 0&1\\1&0\\i&i+k\endmatrix\right)\bigg]<\infty$$
for any $k\geq 1$. (See the proof of Lemma 4.4 of Liggett (1976) for details.)
Next, using the fact that $\nu\in\tilde\Cal I$ implies $\int\tilde\Cal Lg_{m,n}d\nu=0$,
one can show that $\nu$ concentrates on configurations satisfying $g_{m,n}(\eta,\eta')
\leq 1$ for all $m<n$. (See the proofs of Lemmas 4.7 and 4.8 of Liggett (1976).)

The conclusion is the following: In all cases, $\nu$ concentrates on configurations
$(\eta,\eta')$ with the property that there is at most one strict sign
change in the doubly infinite sequence
$$\{...,\eta'(-2)-\eta(-2),\eta'(-1)-\eta(-1),\eta'(0)-\eta(0),\eta'(1)-\eta(1),\eta'(2)-\eta(2),
...\},$$
i.e., such that exactly one of the following is true:

(a) $\eta=\eta'$,

(b) $\eta\leq\eta'$ and $\eta\neq\eta'$,

(c) $\eta\geq\eta'$ and $\eta\neq\eta'$,

(d) there is a $k$ so that $\eta(i)\leq \eta'(i)$ for all $i\leq k$ with infinitely
many strict inequalities, and $\eta(i)\geq \eta'(i)$ for all $i> k$ with infinitely
many strict inequalities,

(e) there is a $k$ so that $\eta(i)\geq \eta'(i)$ for all $i\leq k$ with infinitely
many strict inequalities, and $\eta(i)\leq \eta'(i)$ for all $i> k$ with infinitely
many strict inequalities.

\noindent (The fact that finitely many strict inequalities is excluded in cases
(d) and (e) is a consequence of the fact that the system is in equilibrium; if there
were finitely many, there would be some rate at which they would be annihilated, and
since they cannot be created, this would contradict stationarity.)
Since each of the above sets of configurations is closed for the 
evolution of the coupled process, if $\nu\in\tilde\Cal I_e$, then $\nu$ concentrates
on exactly one of these sets.

Now take $\mu,\mu'\in \Cal I_e$ such that $\phi(\mu)=\phi(\mu')=0$. By
Proposition 2.14 of Chapter VIII of Liggett (1985), there is a $\nu\in\tilde\Cal I_e$
with marginals $\mu$ and $\mu'$ respectively. We will show next that $\nu$ cannot
concentrate on either of the last two sets described above. Suppose, for example,
that it concentrates on the set described in (d). The $k$ appearing there is random,
so we will denote it by $K(\eta,\eta')$. Let
$$\align u_m=&p_{m-1}\mu\{\eta:\eta(m-1)=1,\eta(m)=0\}=q_{m}\mu\{\eta:\eta(m-1)=0,\eta(m)=1\},\\
u_m'=&p_{m-1}\mu'\{\eta:\eta(m-1)=1,\eta(m)=0\}=q_{m}\mu'\{\eta:\eta(m-1)=0,\eta(m)=1\}.\endalign$$
The second equality in each case comes from the fact that the fluxes are zero.
Then
$$\align\frac{u_m-u_m'}{p_{m-1}}=&\mu\{\eta:\eta(m-1)=1,\eta(m)=0\}-
\mu'\{\eta:\eta(m-1)=1,\eta(m)=0\}\\=
&\nu\left(\matrix 1&1\\1&0\\m-1&m\endmatrix\right)+
\nu\left(\matrix 0&0\\1&0\\m-1&m\endmatrix\right)+\nu\left(\matrix
0&1\\1&0\\m-1&m\endmatrix\right)\\&-\nu\left(\matrix 1&0\\1&1\\m-1&m\endmatrix\right)-
\nu\left(\matrix 1&0\\0&0\\m-1&m\endmatrix\right)-\nu\left(\matrix
1&0\\0&1\\m-1&m\endmatrix\right).\endalign$$
The third term on the right above is zero since $\nu$ concentrates on the
set described in (d). The second, fourth and sixth terms are  bounded by
$$\nu\{(\eta,\eta'):K(\eta,\eta')\leq m\}.\tag 2.5$$
For the sixth term, for example, note that $\eta'(m-1)=1, \eta'(m)=0,\eta(m-1)=0,\eta(m)=1$
implies that $K(\eta,\eta')=m-1$. The probability in (2.5)
tends to zero as $m\rightarrow-\infty$. Therefore,
$$p_{m-1}\bigg[\nu\left(\matrix 1&1\\1&0\\m-1&m\endmatrix\right)-
\nu\left(\matrix 1&0\\0&0\\m-1&m\endmatrix\right)\bigg]-(u_m-u_m')\rightarrow 0$$
as $m\rightarrow-\infty$. Similarly
$$q_{m}\bigg[\nu\left(\matrix 1&1\\0&1\\m-1&m\endmatrix\right)-
\nu\left(\matrix 0&1\\0&0\\m-1&m\endmatrix\right)\bigg]-(u_m-u_m')\rightarrow 0$$
as $m\rightarrow-\infty$. Taking differences, we see that
$$p_{m-1}\bigg[\nu\left(\matrix 1&1\\1&0\\m-1&m\endmatrix\right)-
\nu\left(\matrix 1&0\\0&0\\m-1&m\endmatrix\right)\bigg]-
q_{m}\bigg[\nu\left(\matrix 1&1\\0&1\\m-1&m\endmatrix\right)-
\nu\left(\matrix 0&1\\0&0\\m-1&m\endmatrix\right)\bigg]$$
tends to zero as $m\rightarrow-\infty$. Note that this says that the
sum of the first, second, fifth and sixth terms
on the right side of (2.1) tends to zero as $m\rightarrow-\infty$. 
The third, fourth, seventh and eighth terms are bounded by (2.5), so they
tend to zero individually as $m\rightarrow-\infty$. 
Applying a similar argument to the terms in (2.1) that come from the
right boundary of $[m,n]$, we see that the entire right side
of (2.1) tends to zero as $m\rightarrow-\infty, n\rightarrow\infty$. It follows by passing to the
limit in (2.1) that (2.3) holds, and then that (2.4) holds, contradicting the
assumption that $\nu$ concentrates on the set described in (d). A similar
argument shows that (e) cannot hold either. Therefore, either $\mu\leq\mu'$ or
$\mu'\leq \mu$.

To complete the proof of Theorem 1, take $\mu\in \Cal I_e$ and let $\mu'$ be one of the
extremal reversible measures --- $\nu_{\alpha}$ or $\mu_n$ according to whether (1.2)
holds or not. For each such choice of $\mu'$, we now know that either $\mu\leq\mu'$
or $\mu'\leq \mu$. Suppose first that (1.2) holds.  Then there is an $\alpha^*\in [0,1]$
so that
$$\mu\quad\cases\ \leq\nu_{\alpha}\quad\text{if }\alpha_0\geq\alpha^*\\
\ \geq\nu_{\alpha}\quad\text{if }\alpha_0\leq\alpha^*.\endcases$$
Since the $\nu_{\alpha}$'s, together with the pointmasses on the configurations
that are $\equiv 0$ and $\equiv 1$, form a weakly continuous one parameter family,
it follows that $\mu=\nu_{\alpha^*}$ if $\alpha^*\in (0,1)$, and is the
pointmass on the $\equiv 0$ or $\equiv 1$ configuration if $\alpha^*=0$ or 1.

If (1.2) fails, then
it follows that $\mu$ concentrates on $\cup_nA_n\cup \{0,1\}$, where 
0 and 1 represent the identically 0 and identically 1 configurations respectively, so that the
statement of the theorem holds in this case as well. 
\enddemo

\demo{Proof of Proposition 1} We will prove the first statement; the proof of the
second is similar. Suppose that $p_i\geq \frac 12-\epsilon$ for $k\leq i\leq k+n$, and
take $\mu\in \Cal I$. Summing (1.3)
gives
$$\align n\phi(\mu)=&\sum_{i=k}^{k+n-1}\big[p_i\mu\{\eta:\eta(i)=1,\eta(i+1)=0\}-q_{i+1}
\mu\{\eta:\eta(i)=0,\eta(i+1)=1\}\big]\\
\geq&\sum_{i=k}^{k+n-1}\bigg[\bigg(\frac 12-\epsilon\bigg)\mu\{\eta:\eta(i)=1,\eta(i+1)=0\}-
\bigg(\frac 12+\epsilon\bigg)\mu\{\eta:\eta(i)=0,\eta(i+1)=1\}\bigg]\\=&
-2\epsilon\sum_{i=k}^{k+n-1}\mu\{\eta:\eta(i)=1,\eta(i+1)=0\}\\&\quad+
\bigg(\frac 12+\epsilon\bigg)\sum_{i=k}^{k+n-1}\big[\mu\{\eta:\eta(i)=1\}-
\mu\{\eta:\eta(i+1)=1\}\big]\\\geq&-2\epsilon n+\bigg(\frac
12+\epsilon\bigg)\big[\mu\{\eta:\eta(k)=1\}-
\mu\{\eta:\eta(k+n)=1\}\big].
\endalign$$
Since $n$ is arbitrary, it follows that $\phi(\mu)\geq -2\epsilon$. To remove
the extra factor of two, it suffices to note that
$$\sum_{i=k}^{k+n-1}\mu\{\eta:\eta(i)=1,\eta(i+1)=0\}\leq\frac n2$$
for even $n$, since the events occuring in consecutive sumands above are disjoint.
\enddemo

\heading 3. Sufficient conditions for the existence of nonreversible stationary
distributions\endheading

In this section, we prove Theorem 2 and Corollary 2. 
\demo{Proof of Theorem 2} Consider two choices $p_{m-1},..., p_n$ and $p_{m-1}',..., p_n'$
of jump probabilities satisfying $p_i'\geq p_i$ for each $i$. Quantities related to the
process corresponding to the $p_i'$'s will be identified with a prime. We need to
show that $\phi(\mu_{m,n}')\geq \phi(\mu_{m,n})$. In order to do so, we construct a
coupled process as follows: $(X_t,\eta_t,\eta_t',Y_t)$ has state space
$$\gather\bigg\{(x,\eta,\eta',y)\in Z^1\times\{0,1\}^{[m,n]}\times\{0,1\}^{[m,n]}\times Z^1:x+
\sum_{i=m}^n[\eta(i)-\eta'(i)]+y=0\\\text{ and }x+\sum_{i=m}^k[\eta(i)-\eta'(i)]\geq 0\text{ for all }
m-1\leq k\leq n.\bigg\}
\endgather$$
Transitions inside $[m,n]$ correspond to letting particles in the two configurations move
together as much as possible. Transitions across the boundaries $(m-1,m)$ and $(n,n+1)$
follow the same rules, with $X_t$ and $Y_t$ keeping track of the number of discrepancies
that leave or enter $[m,n]$ at the left or right, respectively. To be more explicit, if the $1\ 0$
appearing below are at sites $i,i+1$ respectively (with $m\leq i<n$), then
$$\left(\matrix \eta'\\\eta\endmatrix\right)=\left(\cdots\matrix 1&0\\1&0\endmatrix
\cdots\right)\rightarrow\cases\left(\cdots\matrix 0&1\\1&0\endmatrix
\cdots\right)\quad&\text{at rate }p_i'-p_i\\
\left(\cdots\matrix 0&1\\0&1\endmatrix
\cdots\right)\quad&\text{at rate }p_i,\endcases$$
while if the  $0\ 1$ appearing below are at
sites $i,i+1$ respectively, for example, then
$$\left(\matrix \eta'\\\eta\endmatrix\right)=\left(\cdots\matrix 0&1\\0&1\endmatrix
\cdots\right)\rightarrow\cases\left(\cdots\matrix 0&1\\1&0\endmatrix
\cdots\right)\quad&\text{at rate }q_{i+1}-q_{i+1}'\\
\left(\cdots\matrix 1&0\\1&0\endmatrix
\cdots\right)\quad&\text{at rate }q_{i+1}'.\endcases$$
Note that in both cases, the discrepancy $\left(\matrix 0\\1\endmatrix\right)$ is
produced to the left of the discrepancy $\left(\matrix 1\\0\endmatrix\right)$.
This means that the inequality $x+\sum_{i=m}^k[\eta(i)-\eta'(i)]\geq 0$ is not
violated by these transitions. 

At the left boundary, one has, for example, the following transitions:
$$\bigg(x,\left(\matrix \eta'\\\eta\endmatrix\right)\bigg)=\bigg(x,\left(\matrix 0\\0\endmatrix
\cdots\right)\bigg)\rightarrow\cases\bigg(x+1,\left(\matrix 1\\0\endmatrix
\cdots\right)\bigg)\quad&\text{at rate }p_{m-1}'-p_{m-1}\\
\bigg(x,\left(\matrix 1\\1\endmatrix
\cdots\right)\bigg)\quad&\text{at rate }p_{m-1},\endcases$$
or
$$\bigg(x,\left(\matrix \eta'\\\eta\endmatrix\right)\bigg)=\bigg(x,\left(\matrix 1\\0\endmatrix
\cdots\right)\bigg)\rightarrow\bigg(x-1,\left(\matrix 1\\1\endmatrix
\cdots\right)\bigg)\quad \text{at rate }p_{m-1}.$$
While we have not listed all the possible transitions, hopefully we have listed enough
so that the reader can easily construct the others. 

Now, start this process at $(0,\eta,\eta,0)$, where $\eta$ is any point in $\{0,1\}^{[m,n]}$. The
limiting distribution as
$t\rightarrow\infty$ of $\eta_t$ is $\mu_{m,n}$ while the limiting distribution of $\eta_t'$ is
$\mu_{m,n}'$. For fixed $m\leq k<n$, let
$$\gather
N_t=\text{ (the number of times a particle in $\eta_t$ has crossed from $k$ to $k+1$ by
time $t$)}\\-\text{ (the number of times a particle in $\eta_t$ has crossed from $k+1$ to $k$ by
time $t$),}\endgather$$
with $N_t'$ being defined in an analogous way in terms of the process $\eta_t'$. Then one
can easily check that
$$N_t'-N_t=X_t+\sum_{i=m}^k[\eta_t(i)-\eta_t'(i)],$$
so that $N_t'\geq N_t$ a.s. On the other hand,
$$\frac d{dt}EN_t=p_kP^{\eta}[\eta_t(k)=1,\eta_t(k+1)=0]-q_{k+1}P^{\eta}[\eta_t(k)=0,\eta_t(k+1)=1],$$
so that
$$\phi(\mu_{m,n})=\lim_{t\rightarrow\infty}\frac d{dt}EN_t=\lim_{t\rightarrow\infty}
\frac{EN_t}t.$$
It follows from these observations that $\phi(\mu_{m,n}')\geq\phi(\mu_{m,n})$ as required.
\enddemo
\demo{Proof of Corollary 2} Let $\mu_{m,n}$ be the stationary measure for the process on
$[m,n]$ corresponding to the given $p_i$'s, and $\mu_{m,n}'$ be the one for the process
with $p_i'\equiv\frac 12+\epsilon.$ By Theorem 2, $\phi(\mu_{m,n})\geq\phi(\mu_{m,n}').$
By Theorem 2.9 of Liggett (1977), 
$$\lim\Sb m\rightarrow-\infty\\n\rightarrow+\infty\endSb\mu_{m,n}'=\nu_{1/2},$$
and therefore,
$$\lim\Sb m\rightarrow-\infty\\n\rightarrow+\infty\endSb\phi(\mu_{m,n}')=\frac{\epsilon}2.$$
It follows that any weak limit $\mu$ of $\mu_{m,n}$ as $m\rightarrow-\infty,n\rightarrow
+\infty$ is in $\Cal I$ and satisfies $\phi(\mu)\geq\frac{\epsilon}2$. The measure
$\mu$ is not reversible since all reversible measures have zero flux.
\enddemo

\heading 4. The noninteracting case\endheading

In this section, we consider a system of independent particles, each of which
evolves as a continuous time Markov chain on $Z^1$ with unit exponential holding
times and transition probabilities $p_i$ and $q_i$ from $i$ to $i+1$ and $i-1$
respectively. It has been known since at least the publication of Doob's classic
1953 book that one way to construct stationary distributions for this system is
to let $\{\eta(i), i\in Z^1\}$ be independent Poisson random variables with
$E\eta(i)=\sigma_i$, where $\sigma=\{\sigma_i\}$ is an invariant measure for the
corresponding one-particle motion:
$$\sigma_i=\sigma_{i-1}p_{i-1}+\sigma_{i+1}q_{i+1}.\tag 4.1$$
By Theorem 4.12 of Liggett (1978), these provide all of the extremal stationary
distributions if the one-particle chain is not positive recurrent.

In such a stationary distribution, the flux of particles between $i$ and $i+1$
is
$$\phi=p_i\sigma_i-q_{i+1}\sigma_{i+1},\tag 4.2$$
which is independent of $i$ by (4.1). 

Solving (4.2) recursively leads to
$$\sigma_n=\frac{p_{n-1}\cdots p_m}{q_n\cdots
q_{m+1}}\sigma_m-\frac{\phi}{q_n}\bigg(1+\frac{p_{n-1}}{q_{n-1}}+\cdots+
\frac{p_{n-1}\cdots p_{m+1}}{q_{n-1}\cdots q_{m+1}}\bigg)$$
for $m<n$. It follows that if $\phi>0$, there is a positive solution $\sigma$ to (4.2) if and
only if 
$$1+\frac{q_0}{p_0}+\frac{q_0 q_1}{p_0 p_1}+\cdots<\infty,$$
while if $\phi<0$, there is a positive solution to (4.2) if and only if
$$1+\frac{p_0}{q_0}+\frac{p_0 p_{-1}}{q_0 q_{-1}}+\cdots<\infty.$$
We therefore have the following result:
\proclaim{Theorem 4} Suppose the $p_i$'s are i.i.d. Then there is an extremal stationary
distribution for the independent particle system with positive flux if
$E\log (p_0/q_0)>0$, and one with negative flux if $E\log (q_0/p_0)>0$.
\endproclaim

Not surprisingly, these are exactly the conditions for a random walk in a random environment
to be transient to the right or left respectively -- see Theorem 1.7 in Solomon (1975).
\bigskip
\centerline{\bf References}
\medskip

\ref\by M. Bramson and T. M. Liggett\paper
Exclusion processes in higher dimensions: Stationary measures
and convergence\jour Ann. Probab.\vol 33\pages
2255--2313 \yr 2005\endref

\ref\by M. Bramson, T. M. Liggett and T. Mountford\paper
Characterization of stationary measures for one-dimensional
exclusion processes\jour Ann. Probab.\vol 30\pages 1539--1575
\yr 2002\endref

\ref\by J. L. Doob\book Stochastic Processes\yr 1953\publ Wiley\endref

\ref\by R. J. Harris and R. B. Stinchcombe\paper Disordered asymmetric
simple exclusion process: Mean-field treatment\jour Phys. Rev. E\vol 70 \pages1--15
\yr 2004\endref

\ref\by P. Jung\paper Extremal reversible measures for the
exclusion process\jour J. Stat. Phys.\vol 112\pages 165--191
\yr 2003\endref

\ref \by T. M. Liggett\paper Coupling the simple exclusion process
\jour Ann. Probab.\vol 4\pages 339--356\yr 1976\endref

\ref \by T. M. Liggett\paper Ergodic theorems for the asymmetric simple exclusion
process II
\jour Ann. Probab.\vol 5\pages 795--801\yr 1977\endref

\ref \by T. M. Liggett\paper Random invariant measures for Markov chains,
and independent particle systems
\jour Z. Wahr. verw. Geb.\vol 45\pages 297--313\yr 1978\endref

\ref \by T. M. Liggett\book Interacting Particle Systems\publ Springer-Verlag
\publaddr New York\yr 1985\endref

\ref \by F. Solomon\paper Random walks in a random environment\jour Ann. Probab.
\vol 3\pages 1--31\yr 1975\endref

\ref \by G. Tripathy and M. Barma\paper
Driven lattice gases with quenched disorder: Exact results and  
different macroscopic regime\jour
Phys. Rev. E \vol 58 \pages 1911 --  1926 \yr1998\endref

\bigskip

\noindent Department of Mathematics

\noindent University of California, Los Angeles

\noindent 405 Hilgard Ave.

\noindent Los Angeles CA 90095
\bigskip
\noindent email: lchayes\@math.ucla.edu; tml\@math.ucla.edu

\noindent URL: http://www.math.ucla.edu/\~\ prob-mp/index.html; http://www.math.ucla.edu/\~\ tml/
\end